\input amstex
\documentstyle{amsppt}
\magnification1200 
\tolerance=1000
\def\n#1{\Bbb #1}

\def\lim{\hbox{lim }}

\def\Hom{\hbox{Hom }}

\def\dim{\hbox{dim }}

\def\ve{\varepsilon}
\def\vf{\varphi}

\topmatter
\title The set of jumping conics of a locally free sheaf of dimension 2 on $P^2$
\endtitle 
\author
Dmitry Logachev 
\endauthor
\NoRunningHeads 
\date June 1981 \enddate
\thanks The author thanks his advisor A.S.Tikhomirov for the guidance, and A.N.Tyurin and V.A.Iskovskikh for many valuable remarks. \endthanks 
\abstract We consider a locally free sheaf $F$ of dimension 2 on $P^2$. A conic $q$ on $P^2$ is called a jumping conic if the restriction of $F$ to $q$ is not the generic one. We prove that the set of jumping conics is the maximal determinantal variety of a skew form. Some properties of this skew form are found. Translation from Russian; the original is publshed in: "Constructive algebraic geometry". Yaroslavl, 1981. N.194, p. 79 -- 82. 
 \endabstract 
\address 
Department of Mathematics
Yaroslavl State Pedagogical University 
Yaroslavl, Russia
\endaddress
\endtopmatter
\document 
Let $V$ be a vector space of dimension 3 over $\n C$, $P^2=P(V)$ its projectivization, $F$ a locally free sheaf of dimension 2 on $P^2$ satisfying $c_1(F)=0$,  $c_2(F)=n$, $H^0(P^2, F(1))=0$. Let $q$ be a non-singular conic in $P^2$. We have $F|_q=O(d)\oplus O(-d)$. For a generic $q$ we have $d=0$, let us call these conics regular. If $d>0$  we shall call them jumping conics, the number $d$ is called the multiplicity of jump of $F$ at $q$. 

The main result of the present paper is
\medskip
{\bf Theorem 1.} The set of jumping conics is the maximal determinantal variety of a skew form on $H^1(F\otimes \Omega(-1))$. 

For the reader's convenience, before giving a proof we give firstly the main construction not at all variety of conics on $P^2$ but at a given fixed conic. The solution of the analogous problem in [B] uses a difference between the values of $h^1(F(-1)|_l)$ for an ordinary straight line $l$ and a jumping line. It is easy to see that $\forall k,l$ the numbers $h^i(F(k)|_q)$ coincide for an ordinary conic and a conic of simple jump. There is no sheaf $O(\frac12)$, so we use the sheaf $\Omega(1)$. Let us denote $E_k:=H^1(F\otimes \Omega(k))$. We have $\dim E_{-1}=\dim E_{1}=2n$, and let $t: E_{-1} \otimes S^2V^* \to E_1$ be the $\cup$-multiplication. 
\medskip
{\bf Proposition 2. } $\forall q \in S^2V^*$ the map $t_q:E_{-1}\to E_{1}$ is a skew symmetric bilinear form with respect to a duality $<.,.>: E_{-1}\otimes E_{1}\to \n C$.
\medskip
{\bf Proof.} The duality is defined as follows. Formulas $\lambda^2(F)=O$, $\lambda^2(\Omega)=O(-3)$ imply existence of maps $\ve_1: F\otimes F\to O$, $\ve_2: \Omega\otimes \Omega\to O(-3)$ which are skew symmetric in the following meaning: the diagrams 

$$\matrix  F\otimes F & & \\ & \overset{\ve_1}\to{\searrow} & \\ \downarrow \sigma_1 && O\\ & \overset{\ve_1}\to{\nearrow} & \\  F\otimes F & &\endmatrix  \ \ \  \ \ \ \ \ \ \  \ \ \ \   \matrix \Omega\otimes \Omega & & \\ & \overset{\ve_2}\to{\searrow} & \\ \downarrow \sigma_2 && O\\ & \overset{\ve_2}\to{\nearrow} & \\  \Omega\otimes \Omega & & \endmatrix$$ (here $ \sigma_i$ are transposition maps, i.e. $ \sigma_i(a\otimes b)=b\otimes a$ on local sections) are anti-commutative. Let us construct a map $\ve: F\otimes\Omega\otimes F\otimes\Omega \to O(-3)$ as follows: $\ve(f_1\otimes\omega_1\otimes f_2\otimes\omega_2)= \ve_1(f_1\otimes f_2)\cdot \ve_2(\omega_1\otimes\omega_2)$, where $f_i$, $\omega_i$ are sections. The diagram $$ \matrix  F\otimes \Omega\otimes  F\otimes\Omega & & \\ & \overset{\ve}\to{\searrow} & \\ \downarrow \sigma && O(-3)\\ & \overset{\ve}\to{\nearrow} & \\  F\otimes \Omega\otimes  F\otimes\Omega & & \endmatrix$$ (here $\sigma( f_1\otimes \omega_1\otimes  f_2\otimes\omega_2)= f_2\otimes \omega_2\otimes  f_1\otimes\omega_1$) is commutative. 
\medskip
To complete the proof of Proposition 2, we shall find now $<e_{-1},e_1>=H^2(\ve)(e_{-1}\cup e_1)\in H^2(O(-3))=\n C$. 
\medskip
{\bf Lemma 3.} Skew symmetry of $t_q$ is equivalent to the condition $\forall \  e_{-1}, e'_{-1}\in E_{-1}$ we have $<e_{-1}, t(e'_{-1}\otimes q)>=<e'_{-1}, t(e_{-1}\otimes q)>$.
\medskip
{\bf Proof.} We have $<e_{-1}, t(e'_{-1}\otimes q)>=H^2(\ve)(e_{-1}\cup e'_{-1}\cup q)$ and  

$<e'_{-1}, t(e_{-1}\otimes q)>=H^2(\ve)(e'_{-1}\cup e_{-1}\cup q)$. Cup-product in odd dimensions is anti-commutative, i.e. $H^2(\sigma(-2))(e_{-1}\cup e'_{-1})=-(e'_{-1}\cup e_{-1})$, hence $$H^2(\sigma)(e_{-1}\cup e'_{-1}\cup q)=-(e'_{-1}\cup e_{-1}\cup q) \ \ \ \ \ \  \hbox{and} $$ $$H^2(\ve)(e_{-1}\cup e'_{-1}\cup q)=H^2(\ve)\circ H^2(\sigma)(e_{-1}\cup e'_{-1}\cup q)=-H^2(\ve)(e'_{-1}\cup e_{-1}\cup q)\ \square$$

Let us consider the exact sequence corresponding to the inclusion $q\overset{i}\to{\hookrightarrow}P^2$ 
$$0\to O_{P^2}(-2)\to O_{P^2}\to i_*O_q\to 0 \eqno{(4)}$$
multiply it by $F\otimes\Omega(1)$ and take cohomology: 
$$0\to H^0(F\otimes\Omega(1)|_q)\to E_{-1} \overset{t_q}\to{\to} E_1 \to H^1(F\otimes\Omega(1)|_q)\to 0$$

For regular $q$ (resp. for $q$ of simple jump) we have: $F\otimes\Omega(1)|_q= O(-1)^{\oplus 4}$, resp. $F\otimes\Omega(1)|_q= O(-2)^{\oplus 2}\oplus O^{\oplus 2}$, hence the dimension of both first and fourth terms of (4) are 0, resp. 2. This means that the set of jumping conics is the intersection of $P(S^2(V^*))$ with the maximal determinantal variety of $\Hom_{\hbox{skew}}(E_{-1},E_1)$. Its degree is $n$. 
\medskip
{\bf Proof of Theorem 1.} Let $P^5:=P(S^2(V^*))$ be the set of conics in $P^2$ and $D \hookrightarrow P^5\cdot P^2$ a flag variety defined as follows: $(q,t)\in D \iff t\in q$. We have diagrams

$$\matrix D&\overset{i}\to{\hookrightarrow}&P^5\cdot P^2 \\  \\ &&\pi_5 \swarrow  \searrow\pi_2 \\ \\ && P^5 \ \ \ \ \ \ P^2 \endmatrix \ \ \ \ \ \  \ \ \ \ \ \  \matrix && q_s&\to&s\\  \\ && \downarrow u'_s&&\downarrow u_s\\ \\ P^2& \overset{\pi_2\circ i}\to{\leftarrow}& D& \overset{\pi_5\circ i}\to{\to}& P^5 \endmatrix$$ where $s$ is a point of $P^5$, $q_s$ is the corresponding conic, it is the fibre of $\pi_5\circ i$ at $s$, and $(\pi_2\circ i)\circ u'_s$ is simply the inclusion of $q_s$ in $P^2$. 

The exact sequence corresponding to $D$ is 
$$0\to \pi_5^*O_{P^5}(-1)\otimes  \pi_2^*O_{P^2}(-2)\to O_{P^5\cdot P^2} \to i_*O_D\to 0$$ We multiply it by $ \pi_2^*(F\otimes\Omega(1))$ 
$$0\to \pi_5^*O_{P^5}(-1)\otimes  \pi_2^*(F\otimes\Omega(-1))\to \pi_2^*(F\otimes\Omega(1)) \to  \pi_2^*(F\otimes\Omega(1))|_D \to 0$$
and apply $\pi_{5*}$: 
$$0\to \pi_{5*}(\pi_2^*(F\otimes\Omega(1))|_D) \to E_{-1}\otimes O_{P^5}(-1) \overset{\vf}\to{\to} E_{1}\otimes O_{P^5} \to  $$ $$ \to \pi_{5*1}(\pi_2^*(F\otimes\Omega(1))|_D) \to 0$$
Since the functor of restriction to a fibre is right exact, we get that the support of the sheaf 
$ \pi_{5*1}(\pi_2^*(F\otimes\Omega(1))|_D)$ is the set of jumping conics. The restriction of this sheaf to the set of jumping conics is an analog of the sheaf $\theta(1)$ where $\theta$ is the theta-characteristic sheaf for the case of restriction to straight lines ([T]). Its dimension is 2 at conics of simple jump. The sheaf $\pi_{5*}(\pi_2^*(F\otimes\Omega(1))|_D)$ is obviously 0, and the map $\vf$ comes from the $\cup$-multiplication $t: E_{-1}\otimes S^2V^*\to E_1$. $\square$
\medskip
Let us consider some properties of this map. Its composition with the epimorphism $V^*\otimes V^* \to  S^2V^*$ gives us a map $t': E_{-1}\otimes V^*\otimes V^* \to E_1$ which is the composition of two $\cup$-multiplications $ E_{-1}\otimes V^* \to E_0$, $ E_{0}\otimes V^* \to E_1$. The map $t'': E_{-1}\otimes V^* \to E_1 \otimes V $ --- obtained from $t'$ by moving $V^*$ to the right hand side --- can be factored via $E_0$ and hence has the rank $\dim E_0=2n+2$. By analogy with [T] we can choose a basis $\{e_i\}$ of $E_1$ and a basis $\{v_i\}$ of $V^*$ such that the matrix of the map $t''$ in the basis $\{e_i\otimes v_j\}$ of $E_{-1}\otimes V^*$ is given by the block matriv $A=(A_{ij})$ where $i,j=1,2,3$, $A_{ij}$ is a skew symmetric matrix of size $2n$ and $A_{ij}=A_{ji}$, hence $A$  is skew symmetric. 
\medskip
{\bf References}
\medskip
[B] Barth W., Moduli of vector bundles on the projective line. Inventiones Math., 42 (1977), p. 63 -- 91
\medskip
[T] Tyurin A.N., The geometry of moduli of vector bundles. Russian Math. Surv. 1974, 29:6, p. 57 -- 88. 
\medskip
E-mail: logachev\@{usb.ve}
\enddocument